\documentclass[12pt]{article}
\usepackage[affil-it]{authblk}
\usepackage{amsmath} %Never write a paper without using amsmath for its many new commands
\usepackage{amssymb} %Some extra symbols
\usepackage{graphicx} %If you want to include postscript graphics
\usepackage{enumerate}
\usepackage[utf8]{inputenc}
\usepackage{amsthm}
\usepackage{color}
\usepackage{relsize}

\newtheorem{thm}{Theorem}

\newcommand{\F}{\mathbb{F}}
\newcommand{\Z}{\mathbb{Z}}

\begin{document}

\title{Maximum-order Complexity and Correlation Measures}

\date{}
\author{Leyla I\c s\i k$^1$,  Arne Winterhof$^2$}

\maketitle

\noindent
$^1$ Salzburg University, Hellbrunnerstr. 34,
5020 Salzburg, Austria\\
E-mail: leyla.isik@sbg.ac.at\\

\noindent
$^2$ Johann Radon Institute for Computational and Applied Mathematics\\
Austrian Academy of Sciences, Altenbergerstr.\ 69, 4040 Linz, Austria\\
E-mail: arne.winterhof@oeaw.ac.at

\begin{abstract} 
We estimate the maximum-order complexity of a binary sequence in terms of its correlation measures. 
Roughly speaking, we show that any sequence with small correlation measure up to a sufficiently large order~$k$
cannot have very small maximum-order complexity.
\end{abstract}

\bigskip

{\bf Keywords}:
maximum-order complexity, correlation measure of order $k$, measures of pseudorandomness, cryptography.

\bigskip

{\bf Mathematical Subject Classification}: 11K36, 11T71, 94A55, 94A60.

\section{Introduction}

For a positive integer $N$, the \emph{$N$th linear complexity $L({\cal S},N)$} of a binary sequence ${\cal S}=(s_i)^{\infty}_{i=0}$ is the smallest positive integer 
$L$ such that there are constants $c_0,c_1,...,c_{L-1} \in \F_2$ with
$$s_{i+L}=c_{L-1}s_{i+L-1}+...+c_{0}s_{i},\quad 0\le i \le N-L-1.$$ 
(We use the convention $L({\cal S},N)=0$ if $s_0=\ldots=s_{N-1}=0$ and $L({\cal S},N)=N$ if $s_0=\ldots=s_{N-2}=0\ne s_{N-1}$.)
The $N$th linear complexity is a measure for the predictability of a sequence and thus its unsuitability in cryptography. For surveys on linear complexity and related measures of pseudorandomness see
\cite{gy13,mewi13,ni03,sa07,towi07,wi10}.

Let $k$ be a positive integer. Mauduit and S\'ark\"ozy introduced the ($N$th) \emph{correlation measure of order $k$} of a binary sequence
${\cal S}=(s_i)_{i=0}^\infty$ in \cite{masa97} as
%\begin{equation}\label{legcorr}
$$C_{k}({\cal S},N)=\max_{U,D}\left|\sum^{U-1}_{i=0}(-1)^{s_{i+d_{1}}+s_{i+d_{2}}+...+s_{i+d_{k}}}\right|,$$
%\end{equation}
where the maximum is taken over all $D=(d_1,d_2,...,d_k)$ with non-negative integers $0\leq d_1<d_2<...<d_k$ and $U$ such that $U+d_k\le N$. 
(Actually, \cite{masa97} deals with finite sequences $((-1)^{s_i})_{i=0}^{N-1}$ of length $N$ over $\{-1,+1\}$.)

Brandst\"atter and the second author \cite{branwin06} proved the following relation between the $N$th linear complexity and the correlation measures of order $k$:
\begin{equation}\label{lincorr} L({\cal S},N)\ge N - \max_{1\le k \le L(S,N)+1} C_k({\cal S},N),\quad N\ge 1.
\end{equation}
Roughly speaking, any sequence with small correlation measure up to a sufficiently large order $k$ must have a high $N$th linear complexity as well. 

For example, the Legendre sequence ${\cal L}=(\ell_i)_{i=0}^\infty$ defined by
$$\ell_{i}=\left\{\begin{array}{cl} 1, & \mbox{if $i$ is a quadratic non-residue modulo $p$},\\ 
                                0, & \mbox{otherwise},
                       \end{array}\right.$$
where $p>2$ is a prime, satisfies
$$C_k({\cal L},N)\ll  k p^{1/2} \log p,\quad 1\le N\le p,$$
and thus $(\ref{lincorr})$ implies
$$L({\cal L},N)\gg \frac{\min\{N,p\}}{p^{1/2}\log p},\quad N\ge 1,$$
see \cite{masa97} and \cite[Theorem~9.2]{shbook}.
($f(N)\ll g(N)$ is equivalent to $|f(N)|\le cg(N)$ for some absolute constant $c$.)

The $N$th \emph{maximum-order complexity} $M({\cal S},N)$ of a binary sequence ${\cal S}=(s_i)_{i=0}^\infty$ is the smallest positive integer $M$
such that there is a polynomial $f(x_1,\ldots,x_M)\in \F_2[x_1,\ldots,x_M]$
with
\begin{equation}\label{nonrek} s_{i+M}=f(s_i,s_{i+1},\ldots,s_{i+M-1}),\quad 0\le i\le N-M-1,
\end{equation}
see \cite{ja89,ja91,nixi14}. 
Obviously we have 
$$M({\cal S},N)\le L({\cal S},N)$$
and the maximum-order complexity is a finer measure of pseudorandomness than the linear complexity.

In this paper we analyze the relationship between maximum-order complexity $M({\cal S},N)$ and the correlation measures $C_k({\cal S},N)$ of order $k$. Our main result is the following theorem:
\begin{thm}\label{mainth} For any binary sequence ${\cal S}$ we have
$$M({\cal S},N) \ge N- 2^{M({\cal S},N)+1} \max_{1\le k \le M({\cal S},N)+1} C_k({\cal S},N),\quad N\ge 1.$$
\end{thm}
Again, any nontrivial bound on $C_k({\cal S},N)$ for all $k$ up to a sufficiently large order provides a nontrivial bound on $M({\cal S},N)$.
For example, for the Legendre sequence we get immediately
\begin{equation}\label{legmax} M({\cal S},N)\ge \log(\min\{N,p\}/ p^{1/2}) +O(\log \log p),
\end{equation}
see also \cite[Theorem~9.3]{shbook}.
($f(N)=O(g(N))$ is equivalent to $f(N)\ll g(N)$.) 

We prove Theorem~\ref{mainth} in the next section.

The expected value of the $N$th maximum-order complexity is of order of magnitude $\log N$, see  \cite{ja89} as well as \cite[Remark 4]{nixi14} and references therein.
Moreover, by \cite{al} for a 'random' sequence of length $N$ the correlation measure $C_k({\cal S},N)$ is of order of magnitude $\sqrt{kN\log N}$ and thus by Theorem~\ref{mainth}
$M({\cal S},N)\ge \frac{1}{2} \log N +O(\log\log N)$ which is in good correspondence to the result of \cite{ja89}.

In Section~\ref{further} we mention some straightforward extensions.

\section{Proof of Theorem~\ref{mainth}}

\begin{proof}
Assume ${\cal S}$ satisfies $(\ref{nonrek})$. 
If $s_i=...=s_{i+M-1}=0$ for some $0\le i\le N-M-1$, then $s_{i+M}=f(0,...,0)$. Equivalently,
$(-1)^{s_i}=...=(-1)^{s_{i+M-1}}=1$ implies $(-1)^{s_{i+M}}=(-1)^{f(0,\ldots,0)}$.
Hence, for every $i=0,...,N-M-1$ we have
$$\Big((-1)^{s_{i+M}} - (-1)^{f(0,\ldots,0)}\Big) \prod_{j=0}^{M-1} \Big( (-1)^{s_{i+j}}+1\Big)=0.$$
Summing over $i=0,...,N-M-1$ we get
$$\sum_{i=0}^{N-M-1}\Big((-1)^{s_{i+M}} - (-1)^{f(0,\ldots,0)}\Big) \prod_{j=0}^{M-1} \Big( (-1)^{s_{i+j}}+1\Big)=0.$$
The left-hand side contains one "main" term $\pm(N-M)$ and $2^{M+1}-1$ terms of the form
$$\pm \sum_{i=0} ^{N-M-1} (-1)^{s_{i+j_1}+s_{i+j_2}+\ldots+s_{i+j_k}} $$
with $0\le j_1<j_2<...<j_k\le M$ and $1\le k\le M+1$. Therefore we have
$$N-M\le 2^{M+1} \max_{1 \le k \le M+1} \left |\sum_{i=0} ^{N-M-1} (-1)^{s_{i+j_1}+s_{i+j_2}+\ldots+s_{i+j_k}}\right |$$
and the result follows.
\end{proof}

\section{Further Remarks}
\label{further}

Theorem~\ref{mainth} can be easily extended to $m$-ary sequences with $m>2$ along the lines of \cite{chwi09}:

Let $\xi$ be a primitive $m$th root of unity. Then we have 
$$\sum_{h=0}^{m-1} \xi^{hx}=0\quad \mbox{if and only if}\quad  x\not\equiv 0\bmod m.$$ As in the proof of Theorem~\ref{mainth} we get
$$\sum_{i=0}^{N-M-1} (\xi^{s_{i+M}}-\xi^{f(0,\ldots,0)})\prod_{j=0}^{M-1}\sum_{h=0}^{m-1}\xi^{hs_{i+j}}=0.$$
We have one term of absolute value $N-M$ and $2m^M-1$ terms of the form
\begin{equation}\label{sum} \alpha \sum_{i=0}^{N-M-1} \xi^{h_1s_{i+j_1}+h_2s_{i+j_2}+\ldots+h_ks_{i+j_k}}
\end{equation}
with $1\le h_1,\ldots,h_k<m$, $0\le j_1<j_2<\ldots<j_k\le M$, $1\le k\le M+1$ and $\alpha\in \{1,-\xi^{f(0,\ldots,0)}\}$.

If $m$ is a prime, then  $x\mapsto hx$ is a permutation of $\Z_m$ for any $h\not\equiv 0\bmod m$ and the sums in $(\ref{sum})$ can be estimated by the correlation measure $C_k({\cal S},N)$ of order $k$ 
for $m$-ary sequences as it is defined in \cite{masa02}
and we get 
$$M({\cal S},N)\ge N- 2m^{M({\cal S},N)}\max_{1\le k\le M({\cal S},N)+1} C_k({\cal S},N),\quad N\ge 1.$$
If $m$ is composite, $x\mapsto hx$ is not a permutation of $\Z_m$ if $\gcd(h,m)>1$ and we have to substitute the correlation measure of order $k$ by the power correlation measure of order $k$ introduced
in \cite{chwi09}.\\

Now we return to the case $m=2$.

Even if the correlation measure of order $k$ is large for some small $k$, we may be still able to derive a nontrivial lower bound on the maximum-order complexity by substituting the 
correlation measure of order $k$ by its analog with bounded lags, see \cite{hepawawi15} for the analog of $(\ref{lincorr})$. For example, the two-prime generator ${\cal T}=(t_i)_{i=0}^\infty$, see \cite{branwin05}, 
of length $pq$ with two odd primes $p<q$ satisfies 
$$t_i+t_{i+p}+t_{i+q}+t_{i+p+q}=0$$ if $\gcd(i,pq)=1$ and its correlation measure of order $4$ is obviously close to $pq$, see \cite{ri}.
However, if we bound the lags $d_1<\ldots<d_k<p$ one can derive a nontrivial upper bound on the correlation measure of order $k$  with bounded lags including $k=4$ as well as lower bounds on the maximum-order complexity
using the analog of Theorem~\ref{mainth} with bounded lags.\\

Finally, we mention that the lower bound $(\ref{legmax})$ for the Legendre sequence can be extended to Legendre sequences with polynomials using the results of \cite{go} as well as to their generalization using
squares in arbitrary finite fields (of odd characteristic) using the results of \cite{meya,sawi}. For sequences defined with a character of order $m$ see \cite{masa02}.

\section{Acknowledgement}

The authors are supported by the Austrian Science Fund FWF Projects F5504 and F5511-N26, respectively, which are  part of the Special Research Program "Quasi-Monte Carlo Methods: Theory and Applications". L.I. would like to express her sincere thanks for the hospitality during her visit to RICAM.

\end{document}